\documentclass[12pt]{article}
\usepackage{amsmath, amsfonts, amssymb, amsthm}
\usepackage{verbatim}

\newcommand{\macaulay}{{\sf Macaulay2}}

\theoremstyle{plain}
\newtheorem{theorem}{Theorem}
\newtheorem{proposition}[theorem]{Proposition}

\newtheorem{lemma}[theorem]{Lemma}
\theoremstyle{definition}
\newtheorem{definition}[theorem]{Definition}

\newtheorem{remark}[theorem]{Remark}

\newcommand{\ZZ}{\mathbb{Z}}
\newcommand{\RR}{\mathbb{R}}
\newcommand{\FF}{\mathbb{F}}
\newcommand{\KK}{K}
\DeclareMathOperator{\rank}{rank}
\DeclareMathOperator{\In}{in}
\DeclareMathOperator{\gin}{Gin}
\begin{document}

\date{}
\title{The strong Lefschetz property of the coinvariant ring of the
  %  finite
  Coxeter group of type H$_4$}

\author{
% Toshiaki Maeno (Kyoto University), \\
Yasuhide Numata (Hokkaido University),\\
Akihito Wachi (Hokkaido Institute of Technology)
}

\def\thefootnote{\relax}
\footnotetext{\hspace{-5ex}
}

\maketitle

\begin{abstract}
We prove that the coinvariant ring of the irreducible Coxeter group
of type H$_4$ has the strong Lefschetz property.
\end{abstract}

%#####################################################################
\section{Introduction}
\label{sec:intro}

The strong Lefschetz property (Definition~\ref{defn:slp})
is an abstraction of the hard Lefschetz theorem,
which describes the behavior of the multiplication map 
by the K\"ahler form in the cohomology ring 
of a non-singular algebraic variety.
We show that
the coinvariant ring of the irreducible Coxeter group of type H$_4$
has the strong Lefschetz property (Theorem~\ref{thm:h4-slp}).
This study is a supplement to
\cite{maeno_numata_wachi_slp_of_coxeter},
which determines the set of all strong Lefschetz elements
of the coinvariant rings
of the irreducible Coxeter groups of types other than H$_4$,
in terms of corresponding root systems
(see Remark~\ref{rmk:l-is-in-fundamental-Weyl-chamber}).
For the coinvariant ring of type H$_4$,
Stanley (below Theorem 3.1 of \cite{MR578321})
and Hiller (\cite[Remark on p. 70]{MR630960})
left comments on the difficulty 
in proving the strong Lefschetz property for H$_4$.
To the authors' knowledge it has not been proved up to now.

The difficulty is caused by the complicated structure 
of the Coxeter group of type H$_4$,
and computer algebra systems can not give the answer
under the natural realization of the coinvariant ring
as a quotient ring of $\KK[x_1,x_2,x_3,x_4]$.
The key to our main theorem is to transform the variables 
so that the last variable will be a Lefschetz element.
By this technique, we obtain the theorem 
without heavy computation
except for only a single computation of a Gr\"obner basis.
The computation is executed by
the computer algebra system {\macaulay} \cite{M2}.
The essence of the technique above is paraphrased
as Lemma \ref{lem:criterion-for-slp},
which is a condition for the Lefschetz properties.

Note that
to determine the set of all strong Lefschetz elements
is much more difficult than
to find a strong Lefschetz element,
and we will study the set of the strong Lefschetz elements
for type H$_4$ in a forthcoming paper.

This paper is organized as follows:
In Section~\ref{sec:condition-for-slp} 
we prove Lemma~\ref{lem:criterion-for-slp},
which is the essence of the technique used in our main theorem.
In Section~\ref{sec:H4} we show the main theorem. 
In Section~\ref{sec:remarks} 
we summarize the techniques used in the computation
from the viewpoint of computer algebra systems.

The authors would like to express their sincere thanks 
to Junzo Watanabe for his valuable comments 
during preparation of this paper.

%#####################################################################
\section{A condition for the Lefschetz properties }
\label{sec:condition-for-slp}

In this section we give a necessary and sufficient condition
for graded rings to have the strong Lefschetz property
(Lemma~\ref{lem:criterion-for-slp}).
This lemma is the essence of a technique in proving our main theorem.
The strong or weak Lefschetz property is studied in
\cite{MR951211}, \cite{MR1013668}, 
\cite{MR2033004}\nocite{MR2109963}, \cite{MR1970804}, for instance,
and there are also other conditions for the Lefschetz
properties in terms of generic initial ideals or initial ideals
by Wiebe \cite{MR2111103}.
We see the relation between our lemma and Wiebe's conditions
in the end of this section.

Let $A$ be the polynomial ring $\KK[x_1, x_2, \ldots, x_n]$
over a field $\KK$,
and fix a term order on $A$.
For an ideal $I$ of $A$,
let $\In(I)$ be the {\it initial ideal} of $I$,
which is the monomial ideal generated by the initial monomials
of the polynomials in $I$.
If a monomial in $A$ is not contained in the initial ideal $\In(I)$,
then the monomial is called a {\it standard monomial}
with respect to $I$.
Note that the image of the set of the standard monomials,
under the natural surjection $A \to A/I$,
forms a linear basis of $A/I$.

We recall the Lefschetz properties.
%=====================================================================
\begin{definition}[Watanabe \cite{MR951211}, Iarrobino \cite{MR1184062}]
  \label{defn:slp}
  Let $R$ be a graded ring over a field $\KK$,
  and $R = \bigoplus_{i\ge0} R_i$ its decomposition
  into homogeneous components with $\dim_{\KK} R_i < \infty$.
  The graded ring $R$ is said to have the
  {\it strong (resp.~weak) Lefschetz property},
  if there exists an element $l \in R_1$ such that
  the multiplication map
  $\times l^s: R_i \to R_{i+s}$ ($f \mapsto l^s f$)
  is full-rank for every $i \ge 0$ and $s>0$ (resp.~$s=1$).
  In this case, $l$ is called a {\it Lefschetz element}.

  Suppose that the Hilbert function of the graded ring $R$ is
  {\it symmetric},
  that is, $R = \bigoplus_{i=0}^c R_i$ 
  and $\dim_{\KK} R_i = \dim_{\KK} R_{c-i}$
  for $i = 0, 1, \ldots, \lfloor c/2 \rfloor$.
  In this case,
  it is clear that $R$ has the strong Lefschetz property 
  if and only if there exists $l \in R_1$ and
  $\times l^{c-2i}: R_i \to R_{c-i}$
  is bijective for every $i = 0, 1, \ldots, \lfloor c/2 \rfloor$.
\end{definition}
%=====================================================================

The following lemma is the essence of a technique 
in proving our main theorem.
%=====================================================================
\begin{lemma}
  \label{lem:criterion-for-slp}
  Let $I \subset A = \KK[x_1, x_2, \ldots, x_n]$
  be a homogeneous ideal.
  Take the graded reverse lexicographic order as a term order on $A$.
  Then the following two conditions are equivalent:
  \begin{enumerate}
    \item[(i)]
    The graded ring $A/I$ has 
    the strong (resp.~weak) Lefschetz property,
    and $x_n \bmod I$ is a Lefschetz element.

    \item[(ii)]
    The graded ring $A/\In(I)$ has 
    the strong (resp.~weak) Lefschetz property,
    and $x_n \bmod \In(I)$ is a Lefschetz element.
  \end{enumerate}
\end{lemma}
%=====================================================================
\begin{proof}
When $x_n$ is not necessarily a Lefschetz element,
we claim
\begin{align}
  \label{eq:proof-of-lemma-1}
  \In(I:x_n^s) = \In(I):x_n^s,
\end{align}
for $s > 0$.
It is obvious that $\In(I:x_n^s) \subset \In(I):x_n^s$,
and we prove the other inclusion.
For a monomial $m \in \In(I):x_n^s$,
we can take a homogeneous polynomial $h$ 
such that $m x_n^s + h \in I$,
where $\deg h = \deg (m x_n^s)$ 
and each term of $h$ is smaller than $m x_n^s$.
Thanks to the graded reverse lexicographic order,
$h$ is divisible by $x_n^s$,
and hence $m + h/x_n^s \in I:x_n^s$.
This shows that $\In(I:x_n^s) \supset \In(I):x_n^s$.
We thus have proved (\ref{eq:proof-of-lemma-1}).

When $x_n$ is not necessarily a Lefschetz element,
we have the following formula  
using (\ref{eq:proof-of-lemma-1}).
\begin{align*}
  \rank( \times x_n^s: (A/I)_i \to (A/I)_{i+s} )
  &=
  \dim_{\KK} (A/I)_i - \dim_{\KK} (I:x_n^s/I)_i
  \\ &=
  \dim_{\KK} A_i - \dim_{\KK} (I:x_n^s)_i
  \\ &=
  \dim_{\KK} A_i - \dim_{\KK} (\In(I):x_n^s)_i
  \\ &=
  \rank( \times x_n^s: (A/\In(I))_i \to (A/\In(I))_{i+s} ),
\end{align*}
where the $i$th homogeneous components are denoted as
$(A/I)_i$ and so on.
For a homogeneous ideal $J$ of $A$,
the quotient ring $A/J$ has the strong (resp.~weak) Lefschetz
property with a Lefschetz element $x_n$, 
if and only if the linear map
$\times x_n^s: (A/J)_i \to (A/J)_{i+s}$ is full-rank
for every $i\ge0$ and $s>0$ (resp.~$s=1$).
Therefore it follows from the formula above
that (i) and (ii) of the lemma are equivalent.
\end{proof}
%=====================================================================

In the rest of this section,
We clarify the relation between the lemma and
other conditions for the Lefschetz properties due to Wiebe.
We recall the definition of generic initial ideals.
Fix any term order on $A = \KK[x_1, x_2, \ldots, x_n]$.
For a homogeneous ideal $I$ of $A$,
there exists a Zariski open subset $U \subset GL(n; \KK)$
such that the initial ideals of $\varphi(I)$ are equal to each other
for any $\varphi \in U$.
This initial ideal is uniquely determined,
called the {\it generic initial ideal} of $I$,
and denoted by $\gin(I)$
(see \cite[15.9]{MR1322960}, e.g.).
The following proposition gives other conditions
for the Lefschetz properties.
%=====================================================================
\begin{proposition}[{\cite[Lem 2.7, Prop 2.8, Prop 2.9]{MR2111103}}]
  \label{prop:wiebe}
  Let $I \subset A = \KK[x_1, x_2, \ldots, x_n]$
  be a homogeneous ideal. 
  We have the following:
  \begin{enumerate}
    \item[(i)]
    The graded ring $A/I$ has the strong (resp.~weak) Lefschetz
    property 
    if and only if
    $A/\gin(I)$ has the strong (resp.~weak) Lefschetz property
    with respect to the graded reverse lexicographic order.
    In this case, $x_n \bmod \gin(I)$ is a Lefschetz element
    of $A/\gin(I)$.
    \item[(ii)]
    The graded ring $A/I$ has the strong (resp.~weak) Lefschetz
    property 
    if $A/\In(I)$  has the strong (resp.~weak) Lefschetz property 
    with respect to any term order.
  \end{enumerate}
\end{proposition}
%=====================================================================
The conditions in Proposition~\ref{prop:wiebe}
and the condition in Lemma~\ref{lem:criterion-for-slp} 
relate as follows:
Lemma~\ref{lem:criterion-for-slp} is more practical
than Proposition~\ref{prop:wiebe} (1),
since our lemma does not need generic initial ideals.
Our lemma gives a necessary and sufficient condition
in contrast to Proposition~\ref{prop:wiebe} (2).
In particular, our lemma can be used for checking that an element is
{\it not} a Lefschetz element. 

%#####################################################################
\section{The coinvariant ring of type H$_4$}
\label{sec:H4}
The irreducible Coxeter group of type H$_4$ is of order 14,400
and its root system consists of 120 roots
(see \cite{MR1066460} for details).
In this section,
we show that the coinvariant ring $R$ of
the irreducible  Coxeter group of type H$_4$
has the strong Lefschetz property.
To be concrete, 
we show that an element is a strong Lefschetz element of $R$
by results calculated with the computer algebra system \macaulay.

The {\it coinvariant ring} $R$ has the natural realization
\begin{gather*}
\RR [x_1,x_2,x_3,x_4]/(I_{2k}(x_1,x_2,x_3,x_4) \;|\; k=1,6,10,15 ),
\end{gather*}
where  $I_{2k}(x_1,x_2,x_3,x_4)$ are the polynomials
defined as the sum of the $2k$-th powers of 60 positive roots,
and $I_2, I_{12}, I_{20}, I_{30}$ are the fundamental invariants
of H$_4$ \cite[2.7]{MR926338}.
Note that the root system of type H$_4$ can be realized in
$R_1 = \RR x_1 + \RR x_2 + \RR x_3 + \RR x_4$.
Let
\begin{align*}
\nu_1&= x_1,\\
\nu_2&=  \tau^2 x_1 + x_2 ,\\
\nu_3&=  \tau^4 x_1 + \tau^2 x_2 + x_3 ,\\
\nu_4&=  (\tau^3 + \tau) x_1 + \tau x_2 + x_4 ,
\end{align*}
where $\tau$ is a root of $\tau^2 - \tau - 1$.
We take our candidate $\lambda$ of a strong Lefschetz element of $R$ as
\begin{gather*}
\lambda=\nu_1+\nu_2+\nu_3+\nu_4.   % \label{eq:defofl}
\end{gather*}

\begin{remark}
\label{rmk:l-is-in-fundamental-Weyl-chamber}
When one realizes the root system of type H$_4$ in $R_1$,
the polynomials $\nu_1, \nu_2, \nu_3$ and $\nu_4$ span
the fundamental Weyl chamber with positive coefficients.
In particular, $\lambda$ is in the fundamental Weyl chamber.
We also remark that
\cite{maeno_numata_wachi_slp_of_coxeter} proves that
the set of all strong Lefschetz elements is equal to the union
of all Weyl chambers
for the irreducible Coxeter groups of types other than H$_4$,
and  for type H$_4$ \cite{maeno_numata_wachi_slp_of_coxeter}
only proves that any strong Lefschetz element, if it exists,
is in a Weyl chamber.
\end{remark}

If we use the natural realization of the coinvariant ring as above,
then our computation is too complicated for computer algebra systems
to give the answer.
Thus we transform the ideal
$(I_{2k}(x_1,x_2,x_3,x_4) \;|\; k=1,6,10,15 )$
by the transformation defined by
$\RR[x_1,x_2,x_3,x_4] \to \RR[v_1,v_2,v_3,l]$
$((\nu_1, \nu_2, \nu_3, \lambda) \mapsto (v_1, v_2, v_3, l))$,
and let $I$ be the obtained ideal.
The ideal $I$ is generated by
\begin{align*}
I_{2}& (v_1,v_2 -\tau^2 v_1,v_3 -\tau^2 v_2,l -v_3 -(\tau+1) v_2 -(\tau+1) v_1),\\
I_{12}& (v_1,v_2 -\tau^2 v_1,v_3 -\tau^2 v_2,l -v_3 -(\tau+1) v_2 -(\tau+1) v_1),\\
I_{20}& (v_1,v_2 -\tau^2 v_1,v_3 -\tau^2 v_2,l -v_3 -(\tau+1) v_2 -(\tau+1) v_1),\\
I_{30}& (v_1,v_2 -\tau^2 v_1,v_3 -\tau^2 v_2,l -v_3 -(\tau+1) v_2 -(\tau+1) v_1)
\end{align*}
in the polynomial ring $A = \RR[v_1,v_2,v_3, l]$.
Note that the coinvariant ring $R$
is isomorphic to $A/I$ as graded rings.

We take the reverse lexicographic order such that $v_1>v_2>v_3>l$
as a term order on $A$.
Let $S$ be the set of standard monomials with respect to $I$,
and  $S_{i}$  the set of monomials of degree $i$ in $S$.
By Lemma~\ref{lem:criterion-for-slp} and
the second paragraph of Definition~\ref{defn:slp},
it is enough to show that
\begin{align}
  l^{60-2i}S_i &= S_{60-i} && \text{for all $i<30$.}
  \label{eq:eqtoshow}
\end{align}

We would like to prove Equation~(\ref{eq:eqtoshow}),
but it is difficult to compute Gr\"obner bases
due to intermediate coefficient swells.
To avoid this bottleneck,
we use the finite field $\FF_{13^2}$, i.e.,
the field  obtained
by adjoining a root $\tau$ of $\tau^2-\tau-1$ to  $\FF_{13}$.
%
% Here the characteristic 13 is the minimum prime number $p$
% for which $\tau^2 - \tau - 1$ is irreducible,
% and $\FF_{p^2}[v_1, v_2, v_3, l]/I$
% has the strong Lefschetz property.
%
Namely we prove Equation~(\ref{eq:eqtoshow})
for the ideal $I$ in $\FF_{13^2}[v_1,v_2,v_3,l]$,
and the following \macaulay{} session verifies
Equation~(\ref{eq:eqtoshow}):
%\begin{source}\label{h4prog}
% \verbatiminput{h4slp-m2.txt}
\begin{verbatim}
i1: K =  GF(ZZ/13[tau]/(tau^2-tau-1), Variable => tau);
i2: A = K[v1, v2, v3, l];
i3: v4 = l-v1-v2-v3;
i4: x1 =                            v1;
i5: x2 =                v2 -tau*tau*v1;
i6: x3 =    v3 -tau*tau*v2;
i7: x4 =v4         -tau*v2     -tau*v1;
i8: INVs = { 2*x1, 2*x2, 2*x3, 2*x4,
            x1+x2+x3+x4, x1+x2+x3-x4, x1+x2-x3+x4, x1+x2-x3-x4,
            x1-x2+x3+x4, x1-x2+x3-x4, x1-x2-x3+x4, x1-x2-x3-x4,
            tau*x1+(1/tau)*x2+x3, tau*x1+(1/tau)*x2-x3,
            tau*x1-(1/tau)*x2+x3, tau*x1-(1/tau)*x2-x3,
            tau*x1+(1/tau)*x3+x4, tau*x1+(1/tau)*x3-x4,
            tau*x1-(1/tau)*x3+x4, tau*x1-(1/tau)*x3-x4,
            tau*x1+(1/tau)*x4+x2, tau*x1+(1/tau)*x4-x2,
            tau*x1-(1/tau)*x4+x2, tau*x1-(1/tau)*x4-x2,
            tau*x2+(1/tau)*x4+x3, tau*x2+(1/tau)*x4-x3,
            tau*x2-(1/tau)*x4+x3, tau*x2-(1/tau)*x4-x3,
            (1/tau)*x1+x2+tau*x3, (1/tau)*x1+x2-tau*x3,
            (1/tau)*x1-x2+tau*x3, (1/tau)*x1-x2-tau*x3,
            (1/tau)*x1+x3+tau*x4, (1/tau)*x1+x3-tau*x4,
            (1/tau)*x1-x3+tau*x4, (1/tau)*x1-x3-tau*x4,
            (1/tau)*x1+x4+tau*x2, (1/tau)*x1+x4-tau*x2,
            (1/tau)*x1-x4+tau*x2, (1/tau)*x1-x4-tau*x2,
            (1/tau)*x2+x4+tau*x3, (1/tau)*x2+x4-tau*x3,
            (1/tau)*x2-x4+tau*x3, (1/tau)*x2-x4-tau*x3,
            x1+tau*x2+(1/tau)*x3, x1+tau*x2-(1/tau)*x3,
            x1-tau*x2+(1/tau)*x3, x1-tau*x2-(1/tau)*x3,
            x1+tau*x3+(1/tau)*x4, x1+tau*x3-(1/tau)*x4,
            x1-tau*x3+(1/tau)*x4, x1-tau*x3-(1/tau)*x4,
            x1+tau*x4+(1/tau)*x2, x1+tau*x4-(1/tau)*x2,
            x1-tau*x4+(1/tau)*x2, x1-tau*x4-(1/tau)*x2,
            x2+tau*x4+(1/tau)*x3, x2+tau*x4-(1/tau)*x3,
            x2-tau*x4+(1/tau)*x3, x2-tau*x4-(1/tau)*x3 };
i9: I2 = k ->(sum (set apply(INVs, lf-> lf^(2*k))));
i10: I = ideal (I2(1), I2(6), I2(10), I2(15));
i11: R = A/I;
i12: S = apply(61,k -> first entries(basis({k},R)));
i13: scan(30,i ->( 
       S' = apply(S_i, m -> m*(l^(60-2*i)) );
       << "  l^" << 60-2*i << " S_" << i << " = S_" << 60-i;
       << " is " << ( S' == S_(60-i) ) << endl;
     ));
\end{verbatim}
In this session, we compute the following:
In {\tt i1}, we define the field $K$ obtained by adjoining a root $\tau$ of $\tau^2-\tau-1$  to  $\FF_{13}$.
In {\tt i2}, we define the polynomial ring $A$.
From {\tt i3} to {\tt i10}, we define the ideal $I$.
In {\tt i11},   we define the coinvariant ring $R=A/I$.
In {\tt i12}, we define the set $S_k$ of the standard monomials of degree $k$.
In {\tt i13}, for each $i<30$, we calculate $S' = l^{60-2i} S_i $  and
compare $S'$ with $S_{60-i}$.
%\end{source}

It is easy to see that the strong Lefschetz property for the
coefficient ring $\FF_{13^2}$
yields that for $\ZZ[\tau]$, and then that for $\RR$.
Thus we have our main result.
%=====================================================================
\begin{theorem}
  \label{thm:h4-slp}
  The coinvariant ring of
  the irreducible Coxeter group of type H$_4$
  has the strong Lefschetz property.
\end{theorem}

%#####################################################################
\section{Final remarks}
\label{sec:remarks}

Here we summarize the techniques that we used for the computation of the
main theorem.
As stated in Section~\ref{sec:H4},
the main technique is to take a Lefschetz element 
as the last variable under the graded reverse lexicographic order.
By this technique we do not need heavy computations 
except for a single computation of a Gr\"obner basis.
Otherwise we need many reductions of polynomials of degrees up to 60
and computations of large determinants with rational entries.
In addition,
we use the finite field $\FF_{13^2}$
instead of the rational number field,
which is a common technique in computer algebra systems.
When we use none or one of these techniques,
the computation is too complicated for computers.
The computer algebra system {\macaulay} returns the answer
only when we use both techniques, 
and the computation takes less than 10 seconds in this case.

%#####################################################################
\bibliographystyle{alpha}
% \bibliography{math}

\end{document}